# ON DECIDING STABILITY OF MULTICLASS QUEUEING NETWORKS UNDER BUFFER PRIORITY SCHEDULING POLICIES

By David Gamarnik[1] and Dmitriy Katz

*Massachusetts Institute of Technology*

One of the basic properties of a queueing network is stability. Roughly speaking, it is the property that the total number of jobs in the network remains bounded as a function of time. One of the key questions related to the stability issue is how to determine the exact conditions under which a given queueing network operating under a given scheduling policy remains stable. While there was much initial progress in addressing this question, most of the results obtained were partial at best and so the complete characterization of stable queueing networks is still lacking.

In this paper, we resolve this open problem, albeit in a somewhat unexpected way. We show that characterizing stable queueing networks is an algorithmically undecidable problem for the case of nonpreemptive static buffer priority scheduling policies and deterministic interarrival and service times. Thus, no constructive characterization of stable queueing networks operating under this class of policies is possible. The result is established for queueing networks with finite and infinite buffer sizes and possibly zero service times, although we conjecture that it also holds in the case of models with only infinite buffers and nonzero service times. Our approach extends an earlier related work [*Math. Oper. Res.* **27** (2002) 272–293] and uses the so-called *counter machine* device as a reduction tool.

**1. Introduction.** Queueing networks are ubiquitous tools for modeling a large variety of real-life processes, such as communication and data networks, manufacturing processes, call centers, service networks and many other real-life systems. It is an important task to design and operate queueing networks so that their performance is acceptable. One of the key qualitative performance measures is stability. Roughly speaking, a queueing network is stable if the total expected number of jobs in the network is bounded

Received August 2007; revised January 2009.
[1]Supported by NSF Grant CMMI-0726733.

as a function of time. In a probabilistic framework, which is typically used to formalize the stability question, it means that the underlying queue length process is positive (Harris) recurrent; see [18, 19, 41]. We do not provide a formal definition of this notion here as, throughout the paper, we consider exclusively deterministic queueing networks, for which stability simply means that the total number of jobs in the network remains bounded as a function of time. The details of the model description and formal definitions of stability are delayed until the next section.

The research on stability questions started with the works of Kumar and Seidman [36], Lu and Kumar [37] and Rybko and Stolyar [45], who, for the first time, identified queueing networks and work-conserving scheduling policies leading to instability, even though every processing unit was nominally underloaded. Namely, the condition $\rho_S < 1$ was satisfied by every server $S$. Here, $\rho_S$ is the average utilization in server $S$ which is measured roughly as a ratio of the total arrival rate into this server to the service rate of this server (see the next section). This initiated the search for tight stability conditions. Important advances were obtained in this direction, most notably the development of the fluid model methodology, which significantly simplifies the stability issue by reducing the underlying stochastic problem to a simpler, deterministic continuous-time continuous-state problem; see [19, 48]. It was established that the stability of the fluid model implies stability of the underlying stochastic network [19, 48] and, partially, the converse result holds as well [20, 30, 40, 43], although not always; see [14, 21]. Yet, even characterizing stability of fluid models turned out to be nontrivial [7, 23, 24, 25] and no full characterization is available either. Meanwhile, it was discovered that certain classes of networks and scheduling policies are universally stable. For example, networks with feedforward (acyclic) structure were proven to be stable under an arbitrary work-conserving scheduling policy; see [18, 19, 22]. First Buffer First Serve, Last Buffer First Serve static buffer priority type scheduling policies were shown to stabilize an arbitrary queueing network satisfying a certain topological restriction (a so-called *re-entrant line*); see [25, 35]. A certain simple scheduling policy based on due dates was shown to stabilize an arbitrary network [15]. The First-In, First-Out (FIFO) policy was proven to be stable in networks where service rates within each server are identical—the so-called *Kelly-type networks* [13]. At the same time, some simple static buffer priority policies are not necessarily stable, as was shown in the original works on instability; see [36, 37, 45]. Also, FIFO policy can lead to instability; see [12, 46].

While most of the aforementioned research activity was conducted in the operations research, electrical engineering and mathematics communities, in parallel and independently, the stability problem was investigated by the theoretical computer science community using the *Adversarial Queueing Network* (AQN) model. The motivation there comes from data networks and the models are somewhat different: no probabilistic assumptions are made on either the arrival or service processes. Instead, an adversary is assumed to inject jobs (communication packets) into the network, which is represented as a graph. The links of the graph serve the roles of processing units and the processing times are typically assumed to be equal to one unit of time deterministically. In this setting, the model is defined to be stable if, for every pattern of packet injections, subject to certain load conditions, the total number of packets remains bounded as a function of time. The AQN was introduced by Borodin et al. [9] and further researched by many authors; see [1, 2, 3, 4, 26, 28, 31, 38, 44, 49]. Many results similar to the stochastic networks counterpart were established. It was shown that while AQN corresponding to an acyclic graph is always stable [9], there are AQN and scheduling policies (usually called *protocols*) which are work-conserving (usually called *greedy*) and which lead to instability; see [1, 31]. It was also established that FIFO can lead to instability [1], even with arbitrary small injection rates [6]. The relevance of fluid models to AQN was established in [26]: stability of the fluid model implies stability of AQN. A partial converse result holds, as was also shown in [26]. Yet, despite impressive progress in the area and interesting parallel development to the stochastic counterpart, tight characterization of stable AQN has still not been achieved.

In this paper, we frame the problem of characterizing stable queueing networks as an algorithmic decision problem: given a queueing network and an appropriately defined scheduling policy, determine whether the network is stable. In order to introduce the problem formally, we consider the simplest possible setting: the interarrival times and service times are assumed to take deterministic rational values. Throughout the paper, we focus exclusively on a simple class of scheduling policies, namely the class of nonpreemptive static buffer priority scheduling policies. We assume that buffers have finite or infinite capacity. Jobs which, upon arrival, see a full (finite) buffer are dropped from the network. Also, we assume that some of the service times can take zero value. The assumptions of finite buffers and zero service times are the only important departures from models studied in the stability literature prior to our work. They are adopted for proof tractability, although we conjecture that our main results remain true in the case of infinite buffer/nonzero service times case as well. The details of the model are given in the following section.

Our main result is that stability of a queueing network operating under a static nonpreemptive buffer priority policy is an undecidable property. Thus, no constructive means of characterizing stable queueing networks for this class of policies is possible. This resolves the open problem of providing tight characterization of stable queueing networks for the class of static nonpreemptive buffer priority policies. Our work extends an earlier work [27] by the first author, were the undecidability result was established for the class of so-called *generalized priority* scheduling policies. Later, this work was extended to the problems of computing stationary distributions and large deviations rates [29]. There are important differences between the current work and [27]. The class of generalized priority policies was not considered in the literature prior to [27]. Additionally, generalized priority policies allow idling, whereas most of the work on stability analysis focuses on work-conserving scheduling policies. Also, [27] considered the single-server setting, whereas, here, we consider the network setting. We note that for the class of buffer priority policies (as well as any other work-conserving scheduling policies), the question of stability of a single server model is trivially decidable: one needs to compute the load factor $\rho$. The system is then stable if and only if $\rho < 1$ ($\rho \leq 1$ if all of the interarrival and service times are deterministic).

The concept of undecidability was introduced in the classical works of Alan Turing in the 1930s and it is one of the principal tools for establishing limitations of certain computational problems. The first problems which were established to be undecidable included the Turing halting problem, the post correspondence problem and several related problems [47]. Typically, one establishes undecidability of a given problem by taking a problem which is already known to be undecidable and establishing a reduction from this problem to the given problem of interest. This method is well known in the computer science literature as the *reduction method*. Lately, several problems were proven to be undecidable in the area of control theory; see [5, 16, 17]. In particular, the work of Blondel et al. [5] used a device known as a *counter machine* or *counter automata* as a reduction tool. In the present paper, as in [5], as well as [27], our proof technique is also based on a reduction from a counter machine model, although the construction details are substantially different from those of [27]. We use a well-known Rybko–Stolyar network [45] as a gadget and construct an elaborate queueing network which is able to emulate the dynamics of an arbitrary counter machine. The undecidability result is then a simple consequence of the undecidability of the halting problem for a counter machine, which is a classical result; see [33].

The remainder of the paper is organized as follows. The model description and the main result are provided in the following section. Background material on a counter machine and undecidability is given in Section 3. Section 4 is devoted to constructing a reduction from a counter machine to a queueing network. Section 5 is devoted to the proof of the main result. It begins with a sketch of the proof, followed by the proof details. In the last subsection of this section, we show that while the condition $\rho_S < 1$ is not satisfied by every server in the network we construct, a simple modification achieves this condition. Some concluding thoughts and questions for further research are given in Section 6.

## 2. Model description and the main result.

2.1. *Deterministic multiclass queueing networks and a static buffer priority scheduling policy.* A multiclass queueing network is described as a collection of $J$ service nodes, $S_1, \ldots, S_J$, and $N$ job classes, $1, 2, \ldots, N$. Each node is assumed to be single-server type and can process at most one job at a time. Each class $i$ is associated with a unique buffer, also denoted by $i$, for convenience. The capacity $B_i$ of the buffer $i$ is finite or infinite and denotes the number of jobs which can be stored in the queue of the class $i$, not including the job in service, if any. The queue length corresponding to class $i$ is the number of jobs in buffer $B_i$ plus possibly (at most one) job currently in service and is denoted by $Q_i(t)$. The total queue length $\sum_{i \in S_j} Q_i(t)$ corresponding to the server $S_j$ at time $t$ is denoted by $Q_{S_j}(t)$.

Each class $i$ is associated with an external arrival process $A_i(0, t)$ which denotes the total number of jobs which arrived externally to the buffer $B_i$ during the time interval $[0, t]$. The arrival processes typically considered in the literature are either random renewal processes (in the stochastic queueing networks literature) or adversarial processes (in the computer science literature). Throughout this paper, we adopt the following simple assumption: the intervals between the arrivals of jobs is a *deterministic* class-dependent rational quantity $a_i$ and the initial delay is some rational $b_i$. Thus, the external arrivals corresponding to the class $i$ occur exactly at times $na_i + b, n = 0, 1, \ldots$, and $A_i(0, t) = \lfloor (t-b)/a_i \rfloor + 1$. The external arrival rate is then $\lambda_i \triangleq 1/a_i$. Let $\lambda = (\lambda_i, \ 1 \leq i \leq N)$. Some classes may not have an associated external arrival process, in which case $a_i = \infty$ ($\lambda_i = 0$) and $A_i(0, t) = 0$ for all $t \geq 0$. We will also write $A_i(t) = 1$ if there is an arrival at time $t$ (i.e., $t = a_i n + b_i$ for some $n \in \mathbb{Z}_+$) and $A_i(t) = 0$ otherwise. Each class $i$ is associated with a deterministic service time $0 \leq m_i < \infty$ which takes a nonnegative rational value. The service rate is $\mu_i \triangleq 1/m_i$. We allow service times to take zero value, namely $\mu_i = \infty$. This assumption is a departure from models considered in prior literature and is adopted for proof tractability. We say that at a given time $t$, server $S$ is busy only if, at time $t$, the server is working on a job which requires a nonzero remaining service time. For every collection of classes $V$, the associated workload $W_V(t)$ at time $t$ is the total time required to serve jobs which are presently in the network and which will *eventually* arrive into classes in $V$, in the absence of new arrivals.

The routing of jobs in the network after the service completions is controlled as follows. A zero–one $N \times N$ sub-stochastic matrix $R$ is fixed. Namely, the row sums of this matrix add up to at most unity and the spectral radius of this matrix is strictly less than unity. For every pair of classes $i, l$ such that $R_{i,l} = 1$, every job which completes service in class $i$ at some time $t$ is immediately routed to buffer $B_l$ after the service completion.

If the buffer is not full, that is, $Q_l(t) < B_l$, then the job is added to the end of the queue in the buffer. If the buffer is full, that is, $Q_l(t) = B_l$, then the job is dropped from the network. The special case $B_l = 0$ is interpreted as follows: a job routed to class $l$ is accepted if and only if the server is idle and can begin processing this job immediately. In fact, the network we will build in Section 4 will only have $B_l = 0$ or $B_l = \infty$. If class $i$ is such that $R_{i,l} = 0$ for all $l$, then the jobs in class $i$ after the service completion depart from the network. Since the routing matrix $R$ has spectral radius less than unity, then $R^m = 0$ for some $m$. Namely, every job leaves the network after some finite number of re-routings. The equation $\bar{\lambda} = \lambda + R^T \bar{\lambda}$, also known as the *traffic equation*, then admits a unique solution, explicitly given as $\bar{\lambda} = [I - P^T]^{-1}\lambda$. Here, $R^T$ denotes the transpose of $R$. For every server $S$, the quantity $\rho_S \triangleq \sum_{i \in S} \bar{\lambda}_i / \mu_i$ is defined to be the traffic intensity or load factor in server $S$.

The selection of jobs for processing is controlled using some *scheduling policy*. In the present paper, we exclusively consider the class of *static non-preemptive buffer priority scheduling policies*. Any such policy $\pi$ is described as follows. For each server $S_j$, a permutation $\theta_j$ of the elements of classes belonging to $S_j$ is fixed. At time $t = 0$ and at every time instance $t$ corresponding to a service completion in $S_j$, the server $S_j$ finds the index $i \in S_j$ with the smallest value $\theta_j(i)$ such that $Q_i(t) > 0$, selects the job in the head of this queue and begins working on it. If $\sum_{i \in S_j} Q_i(s) = 0$, then the server idles until the first time that a job appears in one of the classes and starts working on this job. The vector $\theta = (\theta_j), 1 \leq j \leq J$, then completely specifies the scheduling policy $\pi$. In particular, the scheduling policy is nonpreemptive and nonidling: no service is every interrupted and no server idles whenever at least one of the corresponding queues is nonempty. Static buffer priority policies have been studied extensively in the literature; see [8, 10, 11, 23, 25, 34, 37, 42, 45].

A queueing network, described by servers $S_j, 1 \leq j \leq J$, classes $i = 1, 2, \ldots, N$, the routing matrix $R$, interarrival times $a_i$, delays $b_i$ and service times $m_i$ will be denoted by $\mathcal{Q}$ for brevity. The queueing network $\mathcal{Q}$, together with the scheduling policy $\pi$ and the vector of initial queue lengths $(Q_i(0), 1 \leq i \leq N)$, completely determines the queue length dynamics of the network, namely the vector process $Q(s) = (Q_i(s), s \geq 0)$.

DEFINITION 1. A triplet $(\mathcal{Q}, \pi, Q(0))$ is defined to be stable if

$$\sup_{s \geq 0} \sum_{1 \leq i \leq N} Q_i(s) < \infty. \tag{1}$$

A queueing network $\mathcal{Q}$ together with the scheduling policy $\pi$ is defined to be stable if $(\mathcal{Q}, \pi, Q(0))$ is stable for every $Q(0)$.

When all buffers in the network are infinite, the so-called load condition $\rho_S \leq 1$ for all servers $S$ is necessary for stability. The presence of finite buffers may change the situation, as, for example, the model is trivially stable when all of the buffers are finite. Nevertheless, we will see that the load condition is satisfied by all servers in the specific queueing models we construct in this paper, after appropriate modifications described in Section 5.2.

In models with probabilistic settings, $(Q(s), s \geq 0)$ is typically a stochastic process, in which case the queueing network is defined to be stable if the process is so-called *positive Harris recurrent*; see [18, 19, 41]. Under minor additional assumptions, this implies the property $\sup_{s \geq 0} \sum_{1 \leq i \leq N} \mathbb{E}[Q_i(s)] < \infty$. In our deterministic setting, however, this reduces to the simple condition (1). The main goal of the stability research is developing methods for determining stability of a given triplet $(\mathcal{Q}, \pi, Q(0))$ or pair $(\mathcal{Q}, \pi)$. In many interesting special cases, stability of $(\mathcal{Q}, \pi)$ is implied by stability of $(\mathcal{Q}, \pi, Q(0))$ for a given initial state $Q(0)$. For example, in the stochastic setting, this would be the case provided that the underlying Markov chain is irreducible. Due to the deterministic nature of our model, though, this implication does not necessarily hold and it is important to make the distinction.

2.2. *The main result.* The main result of this paper is establishing the undecidability (noncomputability) of the stability property for the class of buffer priority policies $\theta$. Precisely stated, it is as follows.

THEOREM 1. *The property "$(\mathcal{Q}, \theta, Q(0))$ is stable" is undecidable. Namely, no algorithm can exist which, on every input $(\mathcal{Q}, \theta, Q(0))$, outputs YES if the triplet $(\mathcal{Q}, \theta, Q(0))$ is stable and outputs NO otherwise, where $\mathcal{Q}$ is an arbitrary multiclass queueing network, $\theta$ is an arbitrary nonpreemptive static buffer priority scheduling policy and $Q(0)$ is an arbitrary vector of initial queue lengths.*

To prove Theorem 1, we introduce, in Section 3, a device called a counter machine and its stability. Stability of a counter machine is a property closely related to the so-called halting property, which is a classical undecidable property.

**3. Counter machine, the halting problem and undecidability.** A counter machine (see [5, 33]) is a deterministic computing machine which is a simplified version of a Turing Machine—a formal description of an algorithm performing a certain computational task or solving a certain decision problem. In his classical work on the halting problem, Turing showed that certain decision problems simply cannot have a corresponding solving algorithm and are thus undecidable. For a definition of a Turing Machine and the Turing halting problem, see [47]. Since then, many quite natural problems in mathematics and computer science have been found to be undecidable, Hilbert's tenth problem [39] being one of the most notable examples. The famous Church–Turing thesis states that every computable property can be computed by a Turing Machine. Thus, undecidable problems, that is, problems for which a Turing Machine cannot be built, are truly problems not allowing constructive solutions.

More recently, several undecidability results were obtained in the area of control theory, some of them using a counter machine; see Blondel et al. [5]. For a survey of decidability results in the area of control theory, see Blondel and Tsitsiklis [17]. We also use the counter machine device as our reduction tool and, thus, in the next subsection, we provide a detailed description of a counter machine and state relevant undecidability results.

3.1. *Counter machine and the halting problem.* A counter machine is described by two counters $R_1, R_2$ and a finite collection of states $S$. Each counter $R_i$ contains some nonnegative integer $z_i$ in its register. Depending on the current state $s \in S$ and on whether the content of the registers is positive or zero, the counter machine is updated as follows: the current state $s$ is updated to a new state $s' \in S$ and one of the counters has its number in the register incremented by one, decremented by one or no change in the counters occurs.

Formally, a counter machine is a pair $(S, \Gamma)$. $S = \{s_1, s_2, \ldots, s_m\}$ is a finite set of states and $\Gamma$ is configuration update function $\Gamma : S \times \{0,1\}^2 \to S \times \{(-1,0), (0,-1), (0,0), (1,0), (0,1)\}$. A configuration of a counter machine is an arbitrary triplet $(s, z_1, z_2) \in S \times \mathbb{Z}_+^2$. A configuration $(s, z_1, z_2)$ is updated to a configuration $(s', z_1', z_2')$ as follows. Let $1\{\cdot\}$ be the indicator function. Specifically, for every integer $z$, $1\{z\} = 1$ if $z > 0$ and $1\{z\} = 0$ otherwise. Given the current configuration $(s, z_1, z_2)$, suppose, for example, that $\Gamma(s, 1\{z_1\}, 1\{z_2\}) = (s', 1, 0)$. The current state is then changed from $s$ to $s'$, the content of the first counter is incremented by one and the second counter does not change: $z_1' = z_1 + 1, z_2' = z_2$. We will also write $\Gamma : (s, z_1, z_2) \to (s', z_1 + 1, z_2)$ and $\Gamma : s \to s', \Gamma : z_1 \to z_1 + 1, \Gamma : z_2 \to z_2$. Suppose, on the other hand, that $\Gamma(s, 1\{z_1\}, 1\{z_2\}) = (s', (-1, 0))$. The current state then becomes $s'$, $z_1' = z_1 - 1$, $z_2' = z_2$. Similarly, if $\Gamma(s, b) = (s', (0, 1))$ or $\Gamma(s, b) = (s', (0, -1))$, then the new configuration becomes $(s', z_1, z_2 + 1)$ or $(s', z_1, z_2 - 1)$, respectively. If $\Gamma(s, b) = (s', (0, 0))$, then the state is updated to $s'$, but the contents of the counters do not change. It is assumed that the configuration update function $\Gamma$ is consistent, in the sense that it never attempts to decrement a counter which is equal to zero. The present definition of a counter machine can be extended to the one which incorporates more than two counters, but such an extension is not necessary for our purposes.

Given an initial configuration $(s^0, z_1^0, z_2^0) \in S \times \mathbb{Z}_+^2$, the counter machine uniquely determines the subsequent configurations $(s^1, z_1^1, z_2^1), (s^2, z_1^2, z_2^2), \ldots,$ $(s^t, z_1^t, z_2^t), \ldots$. We fix a certain configuration $(s^*, z_1^*, z_2^*)$ and call it the *halting* configuration. If this configuration is reached, then the process halts and no additional updates are executed. The following theorem establishes the undecidability (also called noncomputability) of the halting property.

THEOREM 2. *Given a counter machine $(S, \Gamma)$, initial configuration $(s^0, z_1^0, z_2^0)$ and the halting configuration $(s^*, z_1^*, z_2^*)$, the problem of determining whether the halting configuration is reached in finite time (the halting problem) is undecidable. It remains undecidable even if the initial and the halting configurations are the same, with both counters equal to zero: $s^0 = s^*, z_1^0 = z_2^0 = z_1^* = z_2^* = 0$.*

The first part of this theorem is a classical result and can be found in [32]. The restricted case of $s^0 = s^*, z_i^0 = z_i^*, i = 1, 2,$ can be similarly proven by extending the set of states and the set of transition rules. It is the restricted case of the theorem which will be used in the current paper.

3.2. *Simplified counter machine (SCM), stability and decidability.* We say that a counter machine is *stable* if the value of counters is bounded as time goes to infinity. Namely, $\sup_t z_1^t < \infty$ and $\sup_t z_2^t < \infty$. It is shown in [26] that determining whether a counter machine which started in a given configuration $(s_1, 0, 0)$ is stable is an undecidable problem, by a simple reduction to the halting problem.

DEFINITION 2. A simplified counter machine (SCM) is a counter machine satisfying the following condition: there exist two functions $\alpha: S \times \{0,1\}^2 \to S, \beta: S \to \{-1, 0, 1\}^2$ such that $\Gamma(s, z_1, z_2) = (\alpha(s, 1\{z_1 > 0\}, 1\{z_2 > 0\}), \beta(\alpha(s, 1\{z_1 > 0\}, 1\{z_2 > 0\})))$. In other words, while the new state $s'$ depends on the entire current configuration $(s, z_1, z_2)$, the incrementing or decrementing of counters at the next step depends only on the *new* state $s'$.

It turns out that this restrictive version of a counter machine is still sufficiently general for our purposes.

PROPOSITION 1. *Given a counter machine, an SCM can be constructed such that the SCM is stable if and only if the given counter machine is stable.*

PROOF. We modify the state space $\{s_j\}, 1 \leq j \leq m,$ to $\{s_j^{\text{odd}}\}_{1 \leq j \leq m} \cup \{(s_j^{\text{even}}, b_1, b_2)\}_{1 \leq j \leq m, b_1, b_2 \in \{-1, 0, 1\}}$. The transition rules are defined as follows: $\alpha(s_j^{\text{odd}}, b_1, b_2) = (s_l^{\text{even}}, \Delta_1, \Delta_2)$ if and only if $\Gamma(s_j, b_1, b_2) = (s_l, \Delta_1, \Delta_2)$, and $\beta(s_l^{\text{even}}, \Delta_1, \Delta_2)) = (\Delta_1, \Delta_2)$. Also, $\alpha(s_l^{\text{even}}, \Delta_1, \Delta_2) = s_l^{\text{odd}}$ and $\beta(s_l^{\text{odd}}) = (0,0)$. It is then not hard to observe that each transition $(s_j, z_1, z_2) \to (s_l, z_1', z_2')$ with $b_1 = z_1' - z_1, b_2 = z_2' - z_2$ is emulated by two transitions in the SCM: $(s_j^{\text{odd}}, z_1, z_2) \to ((s_l^{\text{even}}, b_1, b_2), z_1', z_2') \to (s_l^{\text{odd}}, z_1', z_2')$. □

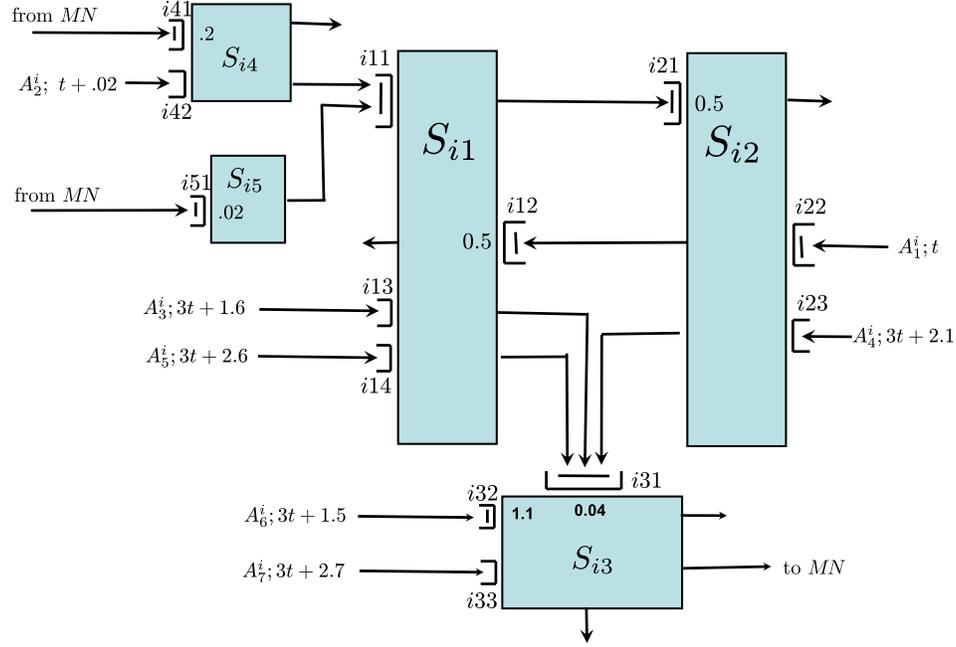

Fig. 1. *Subnetwork* $SN_i$.

COROLLARY 1. *Determining the stability of SCMs with a given initial configuration* $s^*, z_1^* = 0, z_2^* = 0$ *is an undecidable problem.*

**4. Description of the queueing network corresponding to an SCM.** Given an SCM with states $\{s_1, s_2, \ldots, s_m\}$ and counter update rules $\alpha, \beta$, we construct a certain multiclass queueing network, a static buffer priority policy and the vector of queue lengths at time zero. This network, policy and initial state combination will have the property that it is stable if and only if the underlying SCM is stable, thus the reduction goal will be achieved.

We now proceed to the details of the construction. The queueing network consist of three subnetworks denoted, respectively, $SN_1, SN_2$ and $MN$, which stand for *subnetwork* 1, *subnetwork* 2 and the *main network*; see Figures 1 and 2. The subnetwork $SN_i$, $i = 1, 2$, will be in charge of the updates of the counter readings $z_i$. The network $MN$ will be in charge of updating the state $s_i$ of the SCM. We will describe the network structure in detail, as well as the buffer priority scheduling policy implemented in this queueing network. The policy is henceforth denoted by $\theta$. All of the buffer capacities in the network are either zero or infinite.

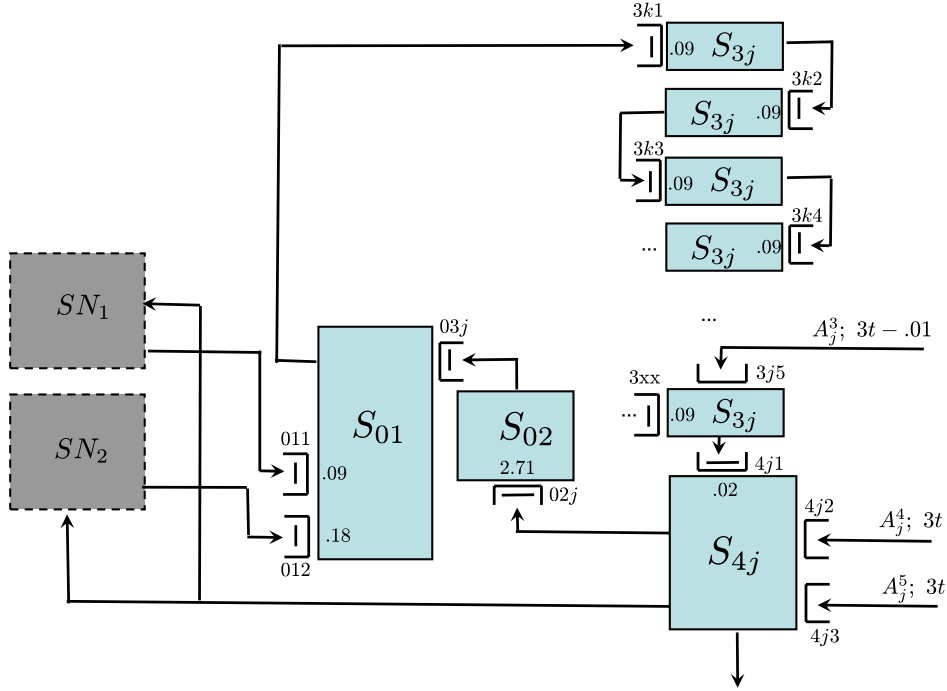

Fig. 2. *Main network MN.*

The subnetworks $SN_i$, $i = 1, 2$, are identical in their topological description. They will only differ in their buffer contents. Hence, we only need to describe one of these subnetworks. In Figures 1 and 2, the buffers with infinite capacity are marked by a vertical bar and the remaining buffers have finite capacity.

4.1. *The description of the subnetwork $SN_i$, $i = 1, 2$.* The subnetwork $SN_i$ consists of five servers, $S_{ij}, j = 1, \ldots, 5$; see Figure 1. The classes (buffers) corresponding to server $S_{ij}$ are denoted by triplets $ijk$. Table 1 lists servers, classes (buffers), the next classes (if any), the corresponding (deterministic) service times, priorities and the buffer capacities. Service times are shown in column 4 and only nonzero service times are shown. If, after service completion, the jobs from a given class exit the system, then the corresponding entry in the next class column is absent. Thus, the unlisted service time entries correspond to zero service time. For each class, we also provide the next class to where the jobs are routed after service completion. If the corresponding entry is empty, it means that the job leaves the network after the service completion. The fifth column corresponds to the priority of this class within the server. For example, the order of priority of classes in server $S_{i1}$ is $i12, i11, i13, i14$, meaning that $i12$ has the highest priority, $i11$ has the next highest priority, etc. The collection of classes $i11, i12, i21, i22$ is defined to be a "Rybko–Stolyar sub-network," or $RSSN_i$. It indeed describes the well-known Rybko–Stolyar network; see [18, 45]. The choice of service times in the subnetwork $SN_i$, as well as in the network $MN$ described in the following section, is somewhat arbitrary, except for service times for classes $i12, i21$ being equal to 0.5. The numbers are arranged so that the proof goes through and is easy to follow. Yet the choice of service times in $i12, i21$ is explained by making the corresponding Rybko–Stolyar network critical, in some appropriate sense. For more details, refer to the beginning of Section 5.

TABLE 1
*Servers and classes in $SN_i$*

| Server | Classes | Next class | Service time | Priority | Capacity |
|---|---|---|---|---|---|
| $S_{i1}$ | $i11$ | $i21$ |  | 2 | $\infty$ |
|  | $i12$ |  | 0.5 | 1 | $\infty$ |
|  | $i13$ | $i31$ |  | 3 | 0 |
|  | $i14$ | $i31$ |  | 4 | 0 |
| $S_{i2}$ | $i21$ |  | 0.5 | 1 | $\infty$ |
|  | $i22$ | $i12$ |  | 2 | $\infty$ |
|  | $i23$ | $i31$ |  | 3 | 0 |
| $S_{i3}$ | $i31$ |  | 0.04 | 2 | $\infty$ |
|  | $i32$ |  | 1.1 | 1 | $\infty$ |
|  | $i33$ | $01i$ of the network $MN$ |  | 3 | 0 |
| $S_{i4}$ | $i41$ |  | 0.2 | 1 | $\infty$ |
|  | $i42$ | $i11$ |  | 2 | 0 |
| $S_{i5}$ | $i51$ | $i11$ | 0.02 | 1 | $\infty$ |

There are seven external arrival processes into subnetwork $SN_i$, denoted by $A^i_j(0, s)$, $j = 1, \ldots, 7$. The corresponding information is summarized in Table 2. For each arrival process, we describe exact arrival times, as well as the class to which the arriving job is routed. For example, the entry $i42$ corresponding to the arrival process $A^i_2$ indicates that jobs arrive precisely at times $0.02, 1.02, 2.02, \ldots$ and are routed to the class $i42$. The arrival times are represented in the form $an + b$ for some explicit constants $a, b$. Here, $a$ is the interarrival time and $b$ is the initial delay. This means that for every nonnegative integer $n$, an arrival occurs at time $an + b$.

4.2. *The description of the main network $MN$.* The main network consists of $2m + 2$ servers, where $m$ is the number of states in the SCM. The servers are $S_{01}$, $S_{02}$, $S_{3j}$, $S_{4j}$, $j = 1, 2, \ldots, m$. The table describing servers, classes, next classes, service times, priorities and buffer capacities is given below as Table 3. The interpretation is the same as for the table for subnetworks $SN_i$. Specific attention is paid to classes $4j3$, $1 \le j \le m$, and the next classes described generically as "$i41, i51$ or exit." The jobs departing from class $4j3$ are routed to:

1. class 141 if $\beta(j) = (-1, 0)$;
2. class 151 if $\beta(j) = (1, 0)$;
3. class 241 if $\beta(j) = (0, -1)$;
4. class 251 if $\beta(j) = (0, 1)$;
5. exit the network if $\beta(j) = (0, 0)$.

TABLE 2
*Arrival processes into $SN_i$*

| Arrival process | Classes | Arrival times |
|---|---|---|
| $A_1^i$ | $i22$ | $n$ |
| $A_2^i$ | $i42$ | $n + 0.02$ |
| $A_3^i$ | $i13$ | $3n + 1.6$ |
| $A_4^i$ | $i23$ | $3n + 2.1$ |
| $A_5^i$ | $i14$ | $3n + 2.6$ |
| $A_6^i$ | $i32$ | $3n + 1.5$ |
| $A_7^i$ | $i33$ | $3n + 2.7$ |

In Table 3, some classes within the same server are assigned the same priority level. This means that the tie is broken arbitrarily. We prefer to assign the same priority level for simplicity. In reality, as we will see, the server will never have to prioritize between these classes as at most one of the corresponding buffers will be nonempty. In order to avoid overcomplicating the figure, the servers $3j$ are described separately for classes $3k1, 3k2, 3k3, 3k4$ and classes $3j5$, although these belong to the same group of servers $3j$, $j = 1, \ldots, m$. Arrivals into the main network are summarized in Table 4. There are $3m$ external arrival processes into subnetwork $MN$, denoted by $A_j^i(0, s)$, $i = 3, 4, 5, j = 1, 2, \ldots, m$. We have started the index $i$ from 3 to avoid confusion with arrival processes $A_j^1, A_j^2$ in networks $SN_i$, $i = 1, 2$. The corresponding information is summarized in Table 4. The arrival times are again represented in the form $an + b$ for some explicit constants $a, b$.

We now describe the initial state of our queueing network at time $s = 0$, namely $Q(0)$. At this time, there is one job in class $02j$ in the main network, where $j$ is such that $s_j = s^*$ is the initial state of the SCM. The service is initiated at time $s = 0$, so the processing of this job will be over at time 2.71. All other buffers in the entire queueing network are empty.

TABLE 3
*Servers and classes in MN*

| Server | Classes | Next classes | Service time | Priority | Capacity |
|---|---|---|---|---|---|
| $S_{01}$ | 011 | | 0.09 | 1 | $\infty$ |
| | 012 | | 0.18 | 2 | $\infty$ |
| | all $03j$ | $3j1$ | | 3 | $\infty$ |
| $S_{02}$ | all $02j$ | $03j$ | 2.71 | 1 | $\infty$ |
| $S_{3j}$ | $3k1$, for all $k$ such that $\alpha(s_k,1,1)=s_j$ | $3k2$ | 0.09 | 1 | $\infty$ |
| | $3k2$, for all $k$ such that $\alpha(s_k,0,1)=s_j$ | $3k3$ | 0.09 | 1 | $\infty$ |
| | $3k3$, for all $k$ such that $\alpha(s_k,1,0)=s_j$ | $3k4$ | 0.09 | 1 | $\infty$ |
| | $3k4$, for all $k$ such that $\alpha(s_k,0,0)=s_j$ | | 0.09 | 1 | $\infty$ |
| | $3j5$ | $4j1$ | 0.02 | 2 | 0 |
| $S_{4j}$ | $4j1$ | | 0.02 | 1 | $\infty$ |
| | $4j2$ | $02j$ | | 2 | 0 |
| | $4j3$ | $i41, i51$ or exit | | 3 | 0 |

**5. Proof of Theorem 1.** Our main result, Theorem 1, follows immediately from Corollary 1 and the following theorem.

THEOREM 3. *The queueing network constructed in the previous section with the prescribed initial state $Q(0)$ is stable if and only if the SCM is stable.*

Before we provide details of the proof of Theorem 3, let us present the overall idea of the proof in the proof sketch below.

PROOF SKETCH OF THEOREM 3. We begin with a brief description of the Rybko–Stolyar network $RSSN_i$, which is embedded in our subnetwork $SN_i, i=1,2$, in relation to servers $S_{i1}, S_{i2}$ and classes $i11, i12, i21, i22$. Instead of two arrival processes feeding class $i11$ in $SN_i$, suppose that we have one external arrival process with arrival times $t=0,1,\dots$. Namely, arrivals occur at the same times as for arrivals into class $i22$. Suppose, as it is in our case, that class $i12$ has priority over class $i11$, and class $i21$ has priority over $i22$. The service times in classes $i11, i12, i21, i22$ are set to take the same

TABLE 4
*Arrival processes into MN*

| Arrival process | Classes | Arrival times |
|---|---|---|
| $A_j^3$ | $3j5$ | $3n-0.01$ |
| $A_j^4$ | $4j2$ | $3n$ |
| $A_j^5$ | $4j3$ | $3n$ |

DECIDING STABILITY 15

values as in our network $SN_i$. Suppose, also, that at time $0^+$, we have $m$ jobs in class $i21$ and no jobs elsewhere. It is a simple exercise to check that at time $m^+$, there will be $m$ jobs in class $i12$ and no jobs elsewhere; at time $(2m)^+$, there will be $m$ jobs in $i21$ and no jobs elsewhere; at time $(3m)^+$, there will be $m$ jobs in $i12$ and no jobs elsewhere, etc. Furthermore, it is a simple exercise to see that the total number of jobs in the four classes $i11, i12, i21, i22$ remains the same $m$ at every integer time $t^+$.

Now, let us go back to our construction. The two Rybko–Stolyar networks $RSSN_i, i = 1, 2$, embedded into $SN_i, i = 1, 2$, will model the two counters in the counter machine, in the sense that the value of the counter $i = 1, 2$ will correspond to roughly the number of jobs in the classes $i11, i12, i21, i22$ at times $3t + 1$ (to be exact, it will correspond to the workload corresponding to these classes; see below). We will arrange the dynamics so that if, during the transition $t \to t + 1$, a counter $i$ has to increment (resp., to decrement, to leave unchanged) its value, then the number of jobs in the Rybko–Stolyar part of $SN_i$ will increase by one (resp., decrease by one, stays the same) over the time period $[3t+1, 3t+4]$. Specifically, say counter $i$ increments its value by one during the transition $t \to t + 1$. We will arrange for exactly one job to arrive from $MN$ into class $i51$ exactly at time $3t + 3$. After an additional delay of 0.02 in server $S_{i5}$, it will arrive into class $i11$ at time $3t + 3.02$. The extra delay of 0.02 is created in order to synchronize with arrivals at time $t + 0.02$ (possibly) coming from class $i42$. The net result is one extra job in the Rybko-Stolyar part of $SN_i$ added during $[3t + 1, 3t + 4]$.

On the other hand, suppose that counter $i$ decrements its value by one during the transition $t \to t + 1$. We will arrange for exactly one job to arrive from $MN$ into $i41$ at time $3t + 3$. This job will occupy server $S_{i4}$ during $(3t + 3, 3t + 3.2)$ and, as a result, the job arriving into class $i42$ at time $3t + 3.02$ will be blocked. The net result (compared to the pure Rybko–Stolyar network described above) is that one job is lost.

The case when the counter does not change simply corresponds to no jobs arriving into $i41$ and $i51$ at time $3t+3$, implying no change in the total number of jobs in the Rybko–Stolyar part of $SN_i$.

Furthermore, the classes $i13, i14, i23$ and classes in the server $S_{i3}$ are constructed so that when a job arrives into zero-buffer class $i33$ at time $3t + 2.7$, it will be processed immediately and sent to $MN$ if the Rybko–Stolyar part of $SN_i$ is empty at time $3t+1$ (namely, counter $i$ is empty) and will be blocked and dropped from the network at time $3t + 2.7$ otherwise. Namely, these classes serve as a testing mechanism for checking whether the counter $i$ is empty or not at time $t$.

Additionally, there is a correspondence between the states of the SCM and the $MN$ network. Specifically, we will arrange that if, at time $t$, the state of SCM is $q$, then, at time $3t$, the server $S_{02}$ will start working on a job in class $02q$. The dynamics is arranged so that if the state of SCM at time $t+1$ is $r$, then, at time $3t+3$, the server $S_{02}$ will start working on a job in class $02r$, thus building the required correspondence between the network $MN$ and the state of the SCM. Specifically, this is arranged as follows. The job in class $02q$ will be processed after 2.71 time units and possibly incur a delay in server $S_{03}$. The delay is either zero, 0.09, 0.18 or 0.27, depending on whether there are jobs arriving into classes 011 and 012 from $SN_1, SN_2$. From the description above, there is a job arriving from $SN_i$ if and only if counter $i$ is empty at time $t$. Thus, the four possible delays uniquely identify which of the counters $i = 1, 2$ are empty and which are not. Next, the job will visit four (possibly repeated) servers among $S_{3j}, 1 \le j \le m$, indexed by four states, $\alpha(q,0,0), \alpha(q,1,0), \alpha(q,0,1), \alpha(q,1,1)$, which can follow state $q$ in the SCM. Depending on the incurred delay, it will be in exactly one of these possible servers at time $3(t+1) - 0.01$ when an external job arrives into this server and is thus blocked. We arrange that it is precisely server $S_{3r}$. The blocked job in buffer $3r5$ is prevented from arriving into class $4r1$ at the same time $3(t+1) - 0.01$ and allows jobs in classes $4r2$ and $4r3$ to be processed at time $3(t+1)$. These will be the only jobs in classes $4j2$ and $4j3$, $j = 1, 2, \ldots, m$, which are served at time $3(t+1)$. One of these jobs arrives into class $02r$, thus completing the cycle and indicating that the new state of the SCM is $r$, and the other job is sent to either $i41$ or $i51$, depending on which of the two counters needs to be updated (if any) and whether the update is increment or decrement. □

For the remainder of the paper, we focus on establishing Theorem 3. We first introduce the following definitions. Let $W_i(s)$ be the combined workload of the servers $S_{i1}, S_{i2}$ in the network $SN_i$ at time $s$. Namely, it is the amount of service required to serve all jobs in servers $S_{i1}, S_{i2}$ at time $s$ when the scheduling policy $\theta$ is implemented. Observe that $W_i(s) = W_{i12}(s) + W_{i21}(s) + 0.5Q_{i22}(s) + 0.5Q_{i11}(s)$, where $W_{i12}(s)$ and $W_{i21}(s)$ stand for the time required to process jobs currently in buffers $i12, i21$ (if any), respectively. We will specifically focus on workloads $W_i(s^-)$, where $s^-$ indicates the time immediately preceding $s$. Thus, if there is an arrival at time $s$, this arrival has not shown up at $s^-$.

For every integer time instance $t = 1, 2, \ldots$, we define the status of the main network $MN$ to be the following quantity: for every $k = 1, 2, \ldots, m$, $Status_{MN}(t) = k$ if, at time $t - 1$, server $S_{02}$ of the network $MN$ started working on a job in class $02k$ and there are no other jobs anywhere in the network $MN$ at time t. Otherwise, $Status_{MN}(t) = -1$.

For each $i = 1, 2$, we also define the status of the subnetwork $SN_i$ at a given time $3t + 1$ for $t \in \mathbb{Z}_+$ as follows. $Status_{SN_i}(3t + 1) = 2W_i((3t + 1)^-)$ if $Q_{i12}(3t + 1)Q_{i21}(3t + 1) = 0$ and there are no jobs anywhere else in the subnetwork $SN_i$, other than possibly in the four classes of $RSSN_i$ (namely, classes $i11, i12, i21, i22$). Otherwise, $Status_{SN_i}(3t+1) = -1$. We do not define $Status_{SN_i}(t)$ at other values of $t$. As we will see shortly, the status functions at time $3t + 1$ will represent the configuration of the SCM at time $t$. Provided that we have initialized our queueing network properly, none of the status functions will ever take value $-1$.

THEOREM 4. *If the configuration of the SCM after $t$ steps is $(s_q, z_1, z_2)$, then $Status_{MN}(3t+1) = q$ and $Status_{SN_i}(3t+1) = z_i$, $i = 1, 2$.*

PROOF. The proof is by induction. For $t = 0$, the statement of Theorem 4 holds because the queueing network initialization makes it so. The remainder of the paper is devoted to proving the induction step. It is given in Section 5.1. □

We now show how this result implies Theorem 3.

PROOF OF THEOREM 3. The idea of the proof is to show that a bound on the value of counters of the SCM implies a bound on the number of jobs in the queueing network at any one time, and vice versa.

Suppose that the SCM is stable. That means that there is a bound $M$ on the maximum value of counters so that $z_1$ and $z_2$ never exceed $M$. Let $(s_j, z_1, z_2)$ be the configuration of the SCM at time $t$. Then, by Theorem 4, at time $(3t+1)^-$, there are $z_1 \leq M$ jobs in $SN_1$, $z_2 \leq M$ jobs in $SN_2$ and one job in the main network. So, at time $(3t+1)^-$, there can be no more than $2M + 1$ jobs in the queueing network. Since there is only a constant number of arrival processes in the network and the arrival process is deterministic, for every time period $[3t + 1, 3(t + 1) + 1)$, the total number of jobs in the network is bounded by $2M + C$ for some constant $C$ which depends only on the network parameters. Thus, if the SCM is stable, so is the queueing network.

Conversely, suppose that the network is stable and that, at any time $t$, the total number of jobs in the network does not exceed $M$ for some finite value $M$. Then, $M$ is also an upper bound on $Status_{SN_i}(3t + 1)$ for every $t$. By Theorem 4, this implies that the values $z_1, z_2$ of the counters of the SCM are bounded by $M$ and therefore the SCM is also stable. □

5.1. *Proof of the induction step of Theorem 4.* This subsection proves the induction step of Theorem 4. Thus, we assume that its statement holds after $t$ steps and prove that it holds after $t + 1$ steps. Assume that the configuration of the SCM at time $t$ is $(s_q, z_1, z_2)$; $Status_{MN}(3t + 1) = q$, $Status_{SN_i}(3t + 1) = z_i$, $i = 1, 2$. Assume that the configuration of SCM at time $t + 1$ is $\Gamma(s_q, z_1, z_2) = (s_r, y_1, y_2)$. We need to show that $Status_{MN}(3t + 4) = r$, $Status_{SN_i}(3t + 4) = y_i$, $i = 1, 2$.

5.1.1. *Dynamics in subnetwork $SN_i$.*

LEMMA 1.   *For every time $s \geq 0$, either $Q_{i12}(s) = 0$ or $Q_{i21}(s) = 0$. Moreover, $\frac{d}{ds}W_i(s) = -1$ whenever $W_i(s) > 0$ and $s \in \mathbb{R}_+$ is not an instance of arrivals into servers $S_{i1}, S_{i2}$.*

REMARK.   The first part of the lemma is a well-known fact from the stability literature, stating that the classes $i12, i21$ constitute a *virtual server* such that only one of the two classes can be served at any given time; see [21, 24].

PROOF OF LEMMA 1.   Suppose that the statement of the lemma does not hold. Then, let $u = \inf(s : Q_{i12}(s) > 0 \text{ and } Q_{i21}(s) > 0)$. That means that both buffers $i12$ and $i21$ are nonempty at time $u^+$, but at least one of the two is empty at time $u^-$. Suppose that this holds for buffer $i12$. This implies that there was an (instantaneous) service completion in buffer $i22$ at time $u$. Class $i21$ has higher priority than class $i22$ (consult Table 1). This implies that the server $S_2$ was not working on the job in class $i21$ at time $u^-$. Since, however, class $i21$ is nonempty at time $u^+$, we conclude that there was an arrival into buffer $i21$ at exactly time $u$. We conclude that there was a simultaneous arrival into buffers $i12$ and $i21$ at time $u$ and buffers $i12$ and $i21$ were empty at time $u^-$.

We now show that such a thing is impossible. Since jobs arrive to $i12$ from $i22$ and into $i22$ from outside at integer times $n$, we see that $u$ must take integer values. We now obtain a contradiction. The jobs arrive into $i11$ only from classes $i42$ and $i51$. Jobs arriving into $i42$ arrive from outside at noninteger times $n + 0.02$. Buffer $i42$ has no capacity and the processing time for this class is zero. Therefore, these jobs can ultimately arrive into $i21$ only at times $n + 0.02$ and not at integer times. Jobs arriving into $i51$ have a nonzero processing time $0.02$. These jobs arrive from the main network $MN$ from classes $4j3$ which correspond to zero capacity buffers and zero processing times. Jobs arrive into $4j3$ from outside at integer times $3n$. Thus, these jobs can ultimately arrive into class $i21$ only at times $3n + 0.02$ and not at integer times. We conclude that jobs cannot ever arrive into $i21$ at integer times.

Similarly, we consider the case where $Q_{i21}(u^-) = 0$. Since $Q_{i21}(u^+) > 0$, there was a service completion in buffer $i11$ at time $u$. We already showed above that this can only occur at times of the form $n + 0.02$. Also, this means that $Q_{12}(u^-) = 0$ since class $i12$ has higher priority than class $i11$. Thus, there was an arrival into $i12$ at time $u$, namely there was a service completion in $i22$ at time $u$. Since $Q_{21}(u^-) = 0$ and the service time in $i22$ is zero, there was an arrival into $i22$ at $u$. But these arrivals only occur at integer times $n$. Again, we obtain a contradiction.

To establish the last part regarding $\frac{d}{ds}W_i(s)$, observe that only jobs in buffers $i12, i21$ have nonzero processing times. Since only one of these buffers can contain a job, the case $W_i(s) > 0$ corresponds to the case of exactly one of these buffers having jobs as, otherwise, if both $i12, i21$ are empty, then the remaining jobs in servers $S_{i1}, S_{i2}$ are processed immediately since they have zero service time requirement. The assertion then follows. □

LEMMA 2. *There are no arrivals into buffers $i41, i51$ during the time interval $[3t+1, 3t+3)$.*

PROOF. Arrivals into classes $i41$ and $i51$ can happen as a result of a departure from one of the classes $4j3$ of the network $MN$. The buffers $4j3$ have zero capacity and zero processing time. Therefore, service completions happen there simultaneously with arrivals from arrival processes $A_j^5$. However, those arrivals occur only at times $3t$. Thus, the first arrival after $3t$ can occur only at time $3t+3$. The assertion then follows. □

LEMMA 3. *During the time interval $[3t+1, 3t+3)$, exactly one of the servers $S_{i1}$ and $S_{i2}$ is busy and $W_i((3t+2)^-) \geq W_i((3t+1)^-)$. In addition, during this time period, jobs in classes $i12$ and $i21$ finish service only at times which are multiples of $0.5$.*

PROOF. By Lemma 1, at most one of servers $S_{i1}, S_{i2}$ does work at any given time. Thus, we need to show that at least one server works during this time period.

By Lemma 2, there are no arrivals into buffers $i41, i51$ during $[3t+1, 3t+3)$. By the inductive assumption, $Status_{SN_i}(3t+1) = z_i \geq 0$, implying, in particular, that there are no jobs in buffer $i41$ at time $3t+1$. Thus, buffer $i41$ is empty during $[3t+1, 3t+3)$. This means that the jobs arriving into class $i42$ at times $3t+1.02$ and $3t+2.02$ will arrive instantly into buffer $i11$. Also, one job will arrive into $i22$ at time $3t+1, 3t+2$. By Lemma 1, only one of the jobs in buffers $i12, i21$ can be served at a time. Thus, the dynamics of the number of jobs in the subnetwork $RSSN_i$ can be viewed as dynamics of a single server queue with service time $0.5$ and arrivals at times $3t+1, 3t+1.02, 3t+2, 3t+2.02$. It is then easy then to explicitly construct $W_i(s)$ during the time period $s \in [3t+1, 3t+3)$, given the initial value $W_i((3t+1)^-)$, and the graph of $W_i(s)$ is depicted in Figures 3–5. The part $[3t+1, 3t+3)$ is identical in all three figures. The differing parts of the graph corresponding to the interval $[3t+3, 3t+4)$ will be used later in Section 5.1.2. In particular, we see that if $W_i((3t+1)^-) > 0$, then $W_i(s)$ is always positive during the time interval $[3t+1, 3t+3)$ and if $W_i((3t+1)^-) = 0$, then $W_i(s)$ is equal to zero only at time $s = 3t+2$. In particular, at least one (and therefore exactly one) of the servers $S_{i1}, S_{i2}$ was busy during the time interval $[3t+1, 3t+3)$. We also see, by inspection, that $W_i((3t+2)^-) \geq W_i((3t+1)^-)$. Finally, by the inductive assumption, $Status_{SN_i}(3t+1) = z_i = 2W_i((3t+1)^-)$; in particular, it is an integer. This means that there is no service in progress in buffers $i12, i21$ at time $3t+1$. Thus, whether or not there are prior jobs in buffers $i12, i21$ at time $3t+1$, there will be service completions exactly at times $3t+1.5, 3t+2, 3t+2.5$ and $3t+3$, as seen by again inspecting Figures 3–5. This proves the second assertion of the lemma. □

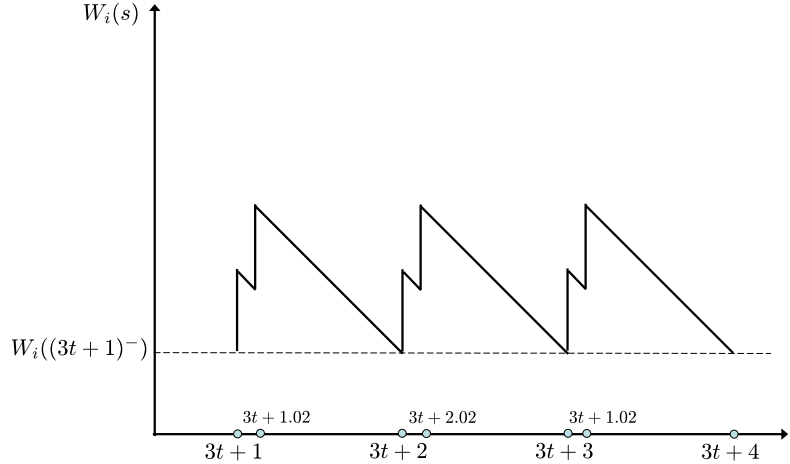

Fig. 3. *Workload $W_i(s)$: case 1.*

LEMMA 4. *Suppose that $Status_{SN_i}(3t+1) \geq 1$. Then, the job $\mathcal{J}$ arriving at time $3t+2.7$ from outside according the arrival process $A_7^i$ will be routed to buffer $01i$ of the network $MN$ at time $3t+2.7$.*

PROOF. At time $3t+1.5$, a job arrives into class $i32$ which requires 1.1 units of processing time. Since $i32$ is the highest priority class in server $S_{i3}$, this server will be busy until time $3t+2.6$. Also, this class having the highest priority implies that there is only one job of this class at a time. Thus, at time $3t+2.6$, buffer $i32$ is empty. Buffer $i31$ has the second highest priority and buffer $i33$, to where the job $\mathcal{J}$ arrives, has the lowest priority. Thus, whether $\mathcal{J}$ will be blocked from service at arrival time $3t+2.7$ depends on the number of jobs in buffer $i31$ at time $3t+2.7$. The processing time for these jobs is 0.04. Therefore, $\mathcal{J}$ will not be blocked if and only if there are at most two jobs in $i31$ since, then, these jobs will be processed not later than $3t+2.6 + 0.04 + 0.04 < 3t+2.7$ and, otherwise, they will be processed at time $3t + 2.6 + 0.04 + 0.04 + 0.04 > 3t + 2.7$. We conclude that $\mathcal{J}$ will be blocked if and only if there are at most two jobs in buffer $i31$. We now show that this is indeed the case provided $Status_{SN_i}(3t+1) \geq 1$.

Jobs arriving into buffer $i31$ depart from classes $i13$, $i14$ and $i23$. These buffers have zero capacity and zero service time. Therefore, they can arrive into $i31$ only at a time of arrival into these three buffers, namely at times $3t + 1.6, 3t + 2.1$ and $3t + 2.6$. In particular, there will be up to three jobs in buffer $i31$ at time $3t + 2.6$. Thus, we need to show that it is impossible for all of these three jobs to arrive into $i31$. We will show that at least one of these jobs is blocked. By Lemma 3, either server $S_{i1}$ or $S_{i2}$ is busy during $[3t+1, 3t+3)$. Suppose that the job arriving into $i13$ at time $3t+1.6$ is not blocked. This means that $S_{i2}$ is busy at time $3t+1.6$. By Lemma 3, it will remain busy until $3t+2$. If it remains busy after this time, then it will remain busy until $3t+2.5$, the job arriving into $i23$ at time $3t+2.1$ is blocked and the assertion is established. Thus, the only remaining possibility is that $S_{i2}$ finishes service at time $3t+2$ and remains idle after this. We will show that a job arriving into $i14$ at time $3t+2.6$ will then be blocked and the proof is then complete. By Lemma 3, $W_i((3t+2)^-) \geq W_i((3t+1)^-) \geq 1$. Thus, there is at least one job in either $S_{i1}$ or $i21$ at time $(3t+2)^-$ which still requires 0.5 units of processing time. We claim that at time $(3t+2)^+$, it is in $i12$. Indeed, it cannot be in $i12$ since the server is idle at this time. For the same reason, it cannot be in $i22$ since service time in this buffer is zero. Also, it cannot be in $i11$ since $S_{i1}$ was idle at $(3t+2)^-$ and the arrivals into $i11$ do not occur at integer times. We conclude that there is at least one job in $i12$ at time $(3t+2)^+$ and no jobs in $i11, i21, i22$ at this time. At time $3t+2$, there is an arrival into $i22$ which then immediately proceeds to $i12$. Thus, we have at least two jobs in $i12$ at time $(3t+2)^+$. The server will work on them during $[3t+2, 3t+3)$ and will block a job arriving into $i14$ at time $3t+2.6$. This completes the proof. $\square$

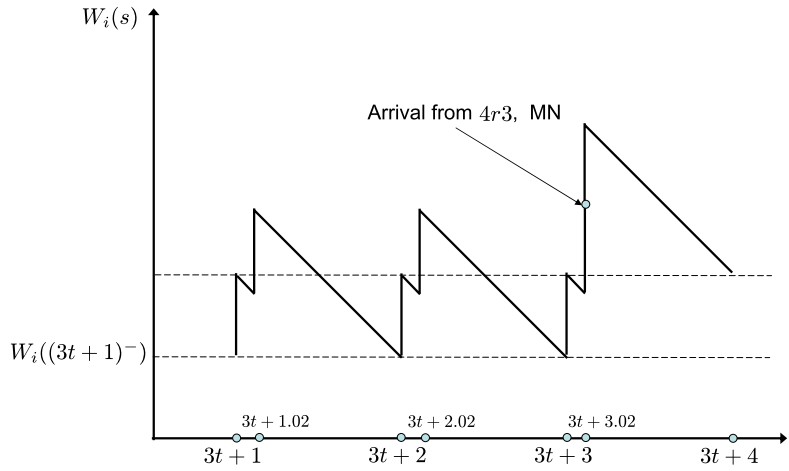

FIG. 4. *Workload $W_i(s)$: case 2.*

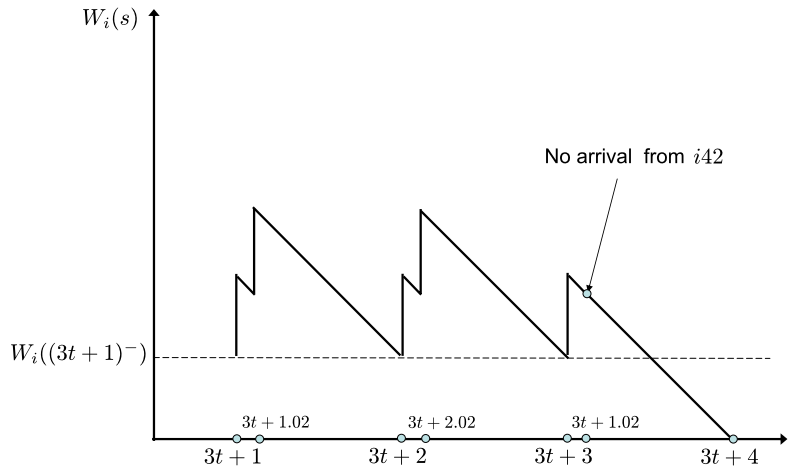

FIG. 5. *Workload $W_i(s)$: case 3.*

LEMMA 5. *Suppose that $Status_{SN_i}(3t+1) = 0$. A job $\mathcal{J}$ arriving at time $3t+2.7$ from outside will then, according to the arrival process $A_7^i$, exit the system immediately.*

PROOF. The proof is very similar to the proof of the previous lemma. We need to show that all three jobs arriving into classes $i13$, $i23$ and $i14$ at times $3t+1.6, 3t+2.1$ and $3t+2.6$, respectively, will not be blocked and will be in buffer $i31$ at time $3t+2.6$. Suppose that $Status_{SN_i}(3t+1) = 0$, that is, $W_i((3t+1)^-) = 0$. The job arriving at time $3t+1$ into buffer $i22$ according to $A_1^i$ will then immediately proceed to buffer $i12$ and occupy server $S_{i1}$ during the time interval $(3t+1, 3t+1.5)$. By Lemma 2, the job arriving into buffer $i24$ at time $3t+1.02$ according to $A_2^i$ will be processed immediately in buffer $i42$ and proceed to buffer $i11$. It will be delayed in buffer $i11$ until $3t+1.5$ and, at this time, will depart to buffer $i21$ and occupy server $S_{i2}$ during the time interval $(3t+1.5, 3t+2)$. Then, again, a job arriving at $3t+2$ into $i22$ will proceed into $i12$ and occupy the server $S_{i1}$ during the time interval $(3t+2, 3t+2.5)$. Finally, the job arriving into $i42$ at time $3t+2.02$ will be delayed in $i11$ until $3t+2.5$ and will then occupy $S_{i2}$ during $(3t+2.5, 3t+3)$. It is clear from this dynamics that all of the three jobs arriving at times $3t+1.6, 3t+2.1$ and $3t+2.6$ into buffers $i13$, $i23$ and $i14$ will be processed immediately and arrive into buffer $i31$ at the same times, $3t+1.6, 3t+2.1$ and $3t+2.6$. $\square$

Combining the results of Lemmas 4 and 5, we obtain the following conclusion.

COROLLARY 2. *If $Status_{SN_i}(3t+1) \geq 1$, then exactly one job arrives into the class $01i$ of network $MN$ at time $3t + 2.7$. If $Status_{SN_i}(3t+1) = 0$, then no job arrives into $01i$ at time $3t + 2.7$.*

5.1.2. *Dynamics in $MN$.* We now switch to the analysis of the dynamics in network $MN$. Recall that, by the inductive assumption $Status_{MN}(3t+1) = q$, we have one job in class $02q$ at time $3t + 1$, which started service at time $3t$, and there are no other jobs in $MN$ at time $3t + 1$. We call this unique job $\mathcal{K}$. Recall that the configuration $(q, x_1, x_2)$ of the SCM at time $t$ is assumed to be updated to the configuration $(r, y_1, y_2)$ at time $t+1$. Introduce $m_1 = \alpha(s_q, 1, 1)$, $m_2 = \alpha(s_q, 0, 1)$, $m_3 = \alpha(s_q, 1, 0)$ and $m_4 = \alpha(s_q, 0, 0)$. Namely, $m_1, m_2, m_3, m_4$ are the four possible values of the state $r$.

LEMMA 6. *During the time interval $(3t + 2.98, 3t + 3.07)$, the job $\mathcal{K}$ will be in server $3r$, buffer $3m_4$ (resp., buffer $3m_3$ or $3m_2$ or $3m_1$) if and only if $x_1 = x_2 = 0$ (resp., if and only if $x_1 = 1, x_2 = 0$ or $x_1 = 0, x_2 = 1$ or $x_1 = x_2 = 0$). This job will leave the network before time $3t + 0.34$.*

PROOF. By the inductive assumption, the job $\mathcal{K}$ will finish service in buffer $02q$ at time $3t + 2.71$ and will arrive into buffer $03q$. It will possibly experience a delay in the corresponding server $S_{01}$ which depends on the presence/absence of jobs in buffers $011, 012$. We now consider four possible cases:

1. Case $x_1 = x_2 = 0$. By the inductive assumption, this means that $Status_{SN_1}(3t+1) = Status_{SN_2}(3t+1) = 0$. By Corollary 2, this means that at time $3t + 2.7$, no jobs arrive into buffers $011, 012$. Since only jobs arriving from buffer $i33$, that is, ultimately from $A_7^i$, can possibly get into buffers $011, 012$, these buffers are empty until at least $3(t+1) + 2.7$. In particular, the job $\mathcal{K}$ arriving into $03q$ at time $3t + 2.71$ will find an idle server and will proceed immediately to buffers $3m_1$, $3m_2$, $3m_3$ and $3m_4$. In each of these buffers, it has the highest priority. Since the service time in each of these buffers is 0.09, it will arrive into these four buffers at exactly the times $3t + 2.71$, $3t + 2.8$, $3t + 2.89$ and $3t + 2.98$. In particular, it will be in buffer $3m_4$ during the time interval $(3t + 2.98, 3t + 3.07)$ and the assertion is established.
2. Case $x_1 = 1, x_2 = 0$. By the inductive assumption, this means that $Status_{SN_1}(3t+1) > 0$, $Status_{SN_2}(3t+1) = 0$. By Corollary 2, this means that at time $3t + 2.7$, no job arrives into buffer $012$ and one job arrives into buffer $011$. This job has the highest priority and requires 0.09 units of processing time. The only difference with the previous case, then, is that the job $\mathcal{K}$ now experiences a delay of 0.09 in server $S_{01}$. Thus, it will arrive into buffers $m_1$, $m_2$, $m_3$ and $m_4$ at exactly the times $3t + 2.8$, $3t + 2.89$, $3t + 2.98$ and $3t + 3.07$. In particular, it will be in the buffer $3m_3$ during the time interval $(3t + 2.98, 3t + 3.07)$ and the assertion is thus established.

3. Case $x_1 = 0, x_2 = 1$. The analysis is similar. We observe that we will have one job in buffer 012 and no jobs in buffer 011 at time $3t + 2.7$. This buffer 012 has the second highest priority; the job $\mathcal{K}$ will experience a delay of 0.18, the processing time of a job in buffer 012.
4. Case $x_1 = x_2 = 1$. The analysis is similar. In this case, we have one job in buffer 011 and one job in buffer 012. The job $\mathcal{K}$ is delayed by $0.18 + 0.09 = 0.27$ time units.

Finally, we again see, by considering the four cases, that the job $\mathcal{K}$ will depart from the network at time $3t + 3.34$, at the latest. This completes the proof of the lemma. $\square$

LEMMA 7. *At time $(3t + 3)^-$, the server $S_{4r}$ is idle and the servers $S_{4j}$, $j \neq r$, are busy processing jobs in buffers $4j1$.*

PROOF. At time $(3t + 3)^-$, the servers $S_{4j}$ can be busy only serving jobs in buffer $4j1$. These jobs arrive from zero capacity buffer $3j5$. These jobs have the highest priority in server $S_{4j}$ and the second highest in $S_{3j}$. Also, these jobs arrive at time $3(t + 1) - 0.01$ into $3j5$. The only way for these jobs to be dropped from zero capacity buffer $3j5$ is by a higher priority buffer in these servers (i.e., one possibly serving job $\mathcal{K}$) being occupied. By Lemma 6, this is the case for exactly one server, namely server $3r$. $\square$

LEMMA 8. $Status_{MN}(3t + 4) = r$.

PROOF. We need to show that at time $3t + 4$, in network $MN$, there is one job in class $02r$ which initiated service at time $3t + 3$ and no jobs elsewhere. By Lemma 6, the job $\mathcal{K}$ will leave the network before time $3t + 3.34 < 3t + 4$. The jobs arriving into zero capacity buffers $4j2, 4j3$, $j \neq r$, at time $3t + 3$ will find, by Lemma 7, a busy server $4j$ and will be dropped from the network. The job arriving into buffer $4r3$ at time $3t + 3$ will find, by Lemma 7, an idle buffer and will immediately proceed to one of the subnetworks $SN_i$. The jobs arriving into buffers $3j5$ at time $3t + 3 - 0.01$ will either be dropped from the network or will proceed to buffers $4j1$ and, after an additional service time 0.02, will leave the network. Thus, they will leave the network before time $3t + 3 + 0.01 < 3t + 4$. We conclude that only the job arriving into buffer $4r2$ at time $3t + 3$ can remain in the network.

By Lemma 7, it will find an idle server $S_{4r}$ and will proceed immediately to buffer $02r$ and begin service there at time $3t + 3$. This completes the proof. □

LEMMA 9. *There are no arrivals into classes $i41$, $i51$ during the time period $[3t+1, 3t+4]$, other than, possibly, at time $3t+3$. At time $3t+3$, at most one job arrives into the four classes $141, 151, 241$ and $251$. Specifically:*

1. $A_{141}(3t+3) = 1$ *if* $\beta(s_r) = (-1, 0)$;
2. $A_{151}(3t+3) = 1$ *if* $\beta(s_r) = (1, 0)$;
3. $A_{241}(3t+3) = 1$ *if* $\beta(s_r) = (0, -1)$;
4. $A_{251}(3t+3) = 1$ *if* $\beta(s_r) = (0, 1)$;
5. *no arrivals if* $\beta(s_r) = (0, 0)$.

PROOF. Arrivals into $i42$ and $i52$ can occur only from buffers $4j3$. These buffers have zero capacity and zero processing times. The arrivals into these buffers occur at times $3n$, $n = 0, 1, \ldots$. By Lemma 7, only server $4r$ will process a job at time $3t + 3$ in buffer $4r3$. According to Table 3 and the corresponding description, it will be routed to one of the buffers $i41$, $i51$ or will leave the network precisely as described by the lemma. □

LEMMA 10. *The following hold for $i = 1, 2$:*

1. $Status_i(3t+4) = Status_i(3t+1)$ *if* $A_{i41}(3t+3) = A_{i51}(3t+3) = 0$;
2. $Status_i(3t+4) = Status_i(3t+1) - 1$ *if* $A_{i41}(3t+3) = 1$;
3. $Status_i(3t+4) = Status_i(3t+1) + 1$ *if* $A_{i51}(3t+3) = 1$.

PROOF. By Lemma 1, we have $Q_{i12}(3t+4)Q_{i21}(3t+4) = 0$. Let us show that at time $3t + 4$, there are no jobs in $SN_i$, other than, possibly, $RSSN_i$. By the inductive assumption, we have $Status_{SN_i}(3t+1) \geq 0$. In particular, at this time, there are no jobs in $SN_i$ outside of $RSSN_i$. We need to show that no jobs arriving during $(3t+1, 3t+4]$ can be outside of $RSSN_i$ at time $3t + 4$.

By Lemma 9, jobs can arrive into $i41, i51$ during $(3t+1, 3t+4]$ only at time $3t + 3$ and only one such job can arrive. Upon arrival, they will experience service time of either 0.2 in $i41$ or 0.02 in buffer $i51$ and they will thus leave the network by time $3t + 3.2$, at the latest.

The jobs arriving into $i42$ at times $3t + 2$, $3t + 3$, $3t + 4$ will either be dropped or proceed to buffer $i11$, which is a part of $RSSN_i$. Thus, at time $3t + 4$, these jobs will either be in $RSSN_i$ or will leave the network (no jobs in $RSSN_i$ feed buffers outside of $RSSN_i$).

We have already analyzed the dynamics of the jobs which arrived into buffers $i13$, $i14$, $i23$, $i32$ and $i33$ at times $3t + 1.6$, $3t + 2.1$ and $3t + 2.6$ as part of the proofs of Lemmas 4 and 5. In particular, we saw that these jobs leave $SN_i$ before time $3t + 2.72$. We have established that there are no jobs in $SN_i$ outside of $RSSN_i$ at time $3t + 4$.

It remains to analyze the value of $Status_i$ at time $3t + 4$. We consider the corresponding three cases:

1. $A_{i41}(3t+3) = A_{i51}(3t+3) = 0$. By Lemma 9, there were no arrivals into classes $i41$, $i51$ in time interval $[3t+1, 3t+4]$. Consider the quantity $W_i(s)$ during this time interval. As long as $W_i(s) > 0$, by Lemma 1, $\frac{d}{ds}W_i(s) = -1$ at time instances $s$ not corresponding to the arrival instances. However, we have arrivals into $i22$ at times $3t + 1$, $3t + 2$ and $3t + 3$, and into $i42$ at times $3t+1+0.02$, $3t+2+0.02$ and $3t+3+0.02$, ensuring that $W_i(s)$ is not 0 for any period of positive length during $[3t+1, 3t+4)$; see Figure 3. In this situation, $W_i(s)$, over the time interval $[3t+1, 3t+4)$, increases by 3 units due to 6 arrivals, and decreases by 3 units due to 6 service completions. Thus, $W_i((3t+4)^-) = W_i((3t+1)^-)$.
2. $A_{i51}(3t+3) = 1$. The job arriving into $i51$ at time $3t + 3$ after a delay of 0.02 will arrive into $i11$, thus increasing $W_i(s)$ by 0.5 at time $s = 3t + 3.02$; see Figure 4. Therefore, $W_i((3t+4)^-) = W_i((3t+1)^-) + 0.5$ and $Status_{SN_i}(3t+4) = Status_{SN_i}(3t+1) + 1$.
3. $A_{i41}(3t+3) = 1$. The job arriving into $i41$ at time $3t + 3$ will occupy server $S_{i4}$ for 0.2 time units. As a result, the job arriving into $i42$ at time $3t + 3.02$ will find a busy server and will be dropped from the network. Comparing this situation with the case $A_{i41}(3t+3) = A_{i51}(3t+3) = 0$ and consulting Figure 5, we obtain the same situation, except that there are no arrivals into $i11$ at time $3t+3.02$. The net result is that $W((3t+4)^-) = W((3t+1)^-) - 0.5$ and $Status_{SN_i}(3t+4) = Status_{SN_i}(3t+1) - 1$.

This completes the proof. □

As an immediate corollary of Lemmas 9 and 10, we obtain the following.

COROLLARY 3. $Status_1(3t + 4) = y_1$ and $Status_2(3t + 4) = y_2$.

Lemma 8 and Corollary 3 prove the induction step for Theorem 4, so its proof is now complete.

5.2. *Load factors.* We will establish below that for some servers in the queueing network constructed in Section 4, the corresponding load factors are greater than unity. As we saw from the proof of our main result, since some of the buffers in our network are finite, overloading some of the servers does not necessarily lead to instability. Yet, this is a significant departure from the standard assumption $\rho_S < 1$ in most of the literature on stability. The goal of this section is to show that simple modifications of our network lead to the same, or a very similar, dynamics, while ensuring the $\rho_S < 1$ condition. Thus, our undecidability result extends to networks with the $\rho_S < 1$ condition satisfied by all servers.

We now compute the load factors $\rho_S$ for each server $S$ encountered in our constructed queueing network and construct appropriate modifications. We begin with the subnetwork $SN_i$. Let us compute the load factors $\rho_{S_{ij}}$, $i = 1, 2$, $j = 1, 2, \ldots, 5$, of the five servers in $SN_i$. The only class in server $S_{i1}$ with nonzero service time (equal to 0.5) is class $i12$. The arrival rate $\bar{\lambda}_{i12}$ into this class equals the external arrival rate into class $i22$, namely $\lambda_{i22} = 1$. Thus, $\rho_{S_{i1}} = 0.5 < 1$ and no modification is needed.

Now, consider server $S_{i2}$. The only class in this server with nonzero service time, equal to 0.5, is class $i21$. The total arrival rate into this class is $\bar{\lambda}_{i21} = \lambda_{i42} + \sum_j \lambda_{4j3}$, where $\lambda_{4j3}$ is the external arrival rate into class $4j3$ in the main network $MN$ and the sum is over all $j$ such that class $4j3$ sends jobs into class $i51$. By construction, $\lambda_{i42} = 1$ and $\lambda_{4j3} = 1/3$. Thus, $\rho_{S_{i2}} \leq (1 + l_1/3)(0.5)$, where $l_1$ is the total number of such classes. As a result, this server is possibly overloaded. We now modify our network as follows. In front of the class $i51$, which is fed by jobs from $MN$, we create a new server with $l_1 + 2$ classes. The first $l_1$ of the classes correspond to arrivals from $MN$ which were originally routed into $i51$. The service rate of these jobs is zero, the buffer size is also zero and, upon service completion, the jobs leave the network. The $(m+1)$st class has external arrivals at exactly the times $3t$ (which are arrival times for classes $4j3$) and service time 0.03. This class has zero buffer and, upon service completion, jobs leave the network. Finally, the class $m + 2$ has arrivals at times $3t + 0.01$, service times 0.01, zero buffer and, upon service completion, jobs are routed into the buffer of the class $i11$. The first $l_1$ classes have the higher priority than class $l_1 + 1$, which, in turn, has higher priority than the class $l_1 + 2$. The load factor of the new server is $(1/3)(0.03 + 0.01) < 1$. Now, let us see how the new server changes the dynamics in the original network. If there is at least one job arriving into classes $1, \ldots, l_1$ in this new server (and we know that only one can arrive at a time), then, since this can only happen at times $3t$, the job arriving into class $l_1 + 1$ is blocked and is dropped from the network. As a result, the job arriving into $l_1 + 2$ at time $3t + 0.01$ is not blocked and is routed into $i11$ at time $3t + 0.02$. On the other hand, if no jobs arrive in classes $1, \ldots, l_1$ at time $3t$, then the job arriving into $l_1 + 1$ at time $3t$ is worked on during the time interval $[3t, 3t + 0.03]$ and blocks the job arriving into $l_1 + 2$ at time $3t + 0.01$, the latter job being dropped from the network. The net effect is the same as when compared with the earlier model: there is one job arriving into $i11$ at time $3t + 0.02$ if and only if there is one job arriving into this class in the original network. But, now, the load factor $\rho_{S_{i2}}$ of the server $S_{i2}$ is $(1 + 1/3)(0.5) < 1$.

Now, consider server $S_{i3}$. We check, in a straightforward way, that $\rho_{S_{i3}} = 3(1/3)(0.04) + (1/3)(1.1) < 1$.

Considering server $S_{i4}$, we see that its load factor, $\rho_{S_{i4}} = l_2(1/3)(0.2)$, may be bigger than unity, where $l_2$ is the total number of classes $4j3$ which may send jobs from $MN$ to class $i41$. Our modification of the network is very simple: replace the service time $0.2$ in $i41$ by $0.2/m$, make arrivals of $A_2^i$ at times $t$, instead of $t + 0.02$, and make service times at $i42$ equal to $0.02$. This makes the load factor of $S_{i4}$ at most $(1/3)m(0.2/m) + 0.02 < 1$. The net effect is the same: if there is an arrival from $MN$ into $i41$, this arrival can occur only at times $3t$ and only one job can arrive at a time. This job occupies the server during $[3t, 3t + 0.2/m]$ and blocks any job arriving into $i42$ according to $A_2^i$ at time $3t$. The latter job is then dropped. If, however, no job arrives into $i41$ at time $3t$, then the job arriving into $i42$ at time $3t$ is processed and, at time $3t + 0.02$, it reaches $i11$, as in the original network.

For server $S_{i5}$, our earlier modification, namely a new server in front of class $i51$, implies that the new load factor is only $\rho_{S_{i5}} = (1/3)(0.02) < 1$.

We now turn to the main network $MN$. Let us compute the load factors $\rho_{S_{01}}, \rho_{S_{02}}, \rho_{S_{3j}}, \rho_{S_{4j}}, 1 \leq j \leq m$, of the servers in $MN$. We have $\rho_{S_{01}} = (1/3)(0.09) + (1/3)(0.18) < 1$ (the two arrival rates $1/3$ are for jobs arriving from subnetworks $SN_1, SN_2$, corresponding to classes $133, 233$). As for server $S_{02}$, we have $\rho_{S_{02}} = m(1/3)(2.71)$ and this server is possibly overloaded as well. We simply replace this server with $m$ identical servers, each dedicated to serving class $02j$, $j = 1, \ldots, m$. Recall that the only function of the server $S_{02}$ was to introduce a fixed delay of $2.71$. Each one of the new $m$ servers has load factor $(1/3)(2.71) < 1$.

Now, let us consider servers $S_{3j}$. We have $\rho_{S_{3j}} = l_3(1/3)(0.09)$, where $l_3$ is the number of classes $03j$ in server $S_{01}$ which can send jobs into server $S_{3j}$. Note that this is also the number of states which can transition into the state $j$ in the SCM. Note that $l_3$ can be as large as $4m$. Thus, this server can be overloaded. Our modification is as follows. Instead of each server $S_{3j}$, we create $4m$ servers $S_{3js}$, $s = 1, \ldots, 4m$. Jobs arriving into classes $3k1$ in server $S_{3j}$ in the original network instead go through servers $S_{3j1}, \ldots, S_{3j(4m)}$, in this order, with service requirement $0.09/(4m)$ in each server. Jobs arriving according to $A_j^3$ into class $3j5$ in the modified version have to go through all of the $4m$ servers $S_{3j1}, \ldots, S_{3j(4m)}$, with zero service time requirement and zero buffer, and are ultimately routed into class $4j1$, as was the case in the original network. It is easy to see that we obtain the same net effect: processing one job in class $3k1$ for $0.09$ time units is replaced by $4m$ subsequent processing stages, each with processing time $0.09/(4m)$. The load factor in each new server is at most $(4m)(0.09)/(4m) < 1$.

Finally, observe that $\rho_{S_{4j}} = (1/3)(0.02) < 1$.

This completes the description of the modified network in which the condition $\rho_S < 1$ is satisfied by every server $S$.

**6. Conclusion.** We have established that there does not exist an algorithm for determining stability of a multiclass queueing network operating under a static nonpreemptive buffer priority scheduling policy. Namely, the underlying problem is undecidable. There are, however, special cases for which the stability can be determined. Characterization of those special cases is of interest. Also of interest is whether our undecidability result holds for FIFO scheduling policy, another frequently studied scheduling policy. Our model incorporated several simplifying assumptions which depart from standard assumptions in the literature on stability of queueing networks. Specifically, we considered networks with possibly finite buffers and zero service times. We have little doubt that the stability property remains undecidable, even for multiclass queueing networks, without these assumptions, but, at present, we do not have a proof of this.

DECIDING STABILITY 31[34] KUMAR, S. and KUMAR, P. R. (1994). Performance bounds for queueing networks and scheduling policies. *IEEE Trans. Automat. Control* **39** 1600–1611. MR1287267

[35] KUMAR, S. and KUMAR, P. R. (2001). Queueing network models in the design and analysis of semiconductor wafer fabs. *IEEE Trans. Robot. Automat.* **17** 548–561.

[36] KUMAR, P. R. and SEIDMAN, T. I. (1990). Dynamic instabilities and stabilization methods in distributed real-time scheduling of manufacturing systems. *IEEE Trans. Automat. Control* **35** 289–298. MR1044023

[37] LU, S. H. and KUMAR, P. R. (1991). Distributed scheduling based on due dates and buffer priorities. *IEEE Trans. Automat. Control* **36** 1406–1416.

[38] LOTKER, Z., PATT-SHAMIR, B. and ROSÉN, A. (2004). New stability results for adversarial queuing. *SIAM J. Comput.* **33** 286–303 (electronic). MR2048442

[39] MATIYASEVICH, Y. (1993). *Hilbert's Tenth Problem*. Nauka, Moscow.

[40] MEYN, S. P. (1995). Transience of multiclass queueing networks via fluid limit models. *Ann. Appl. Probab.* **5** 946–957. MR1384361

[41] MEYN, S. P. and TWEEDIE, R. L. (1993). *Markov Chains and Stochastic Stability*. Springer, London. MR1287609

[42] MORRISON, J. R. and KUMAR, P. R. (1999). New linear program performance bounds for queueing networks. *J. Optim. Theory Appl.* **100** 575–597. MR1684537

[43] PUKHALSKI, A. A. and RYBKO, A. N. (2000). Nonergodicity of queueing networks when their fluid models are unstable. *Problemy Peredachi Informatsii* **36** 26–46. MR1746007

[44] ROSEN, A. (2002). A note on models for non-probabilistic analysis of packet switching networks. *Inform. Process. Lett.* **84** 237–240. MR1931726

[45] RYBKO, A. N. and STOLYAR, A. L. (1992). On the ergodicity of random processes that describe the functioning of open queueing networks. *Problemy Peredachi Informatsii* **28** 3–26. MR1189331

[46] SEIDMAN, T. I. (1994). "First come, first served" can be unstable! *IEEE Trans. Automat. Control* **39** 2166–2171. MR1295752

[47] SIPSER, M. (1997). *Introduction to the Theory of Computability*. PWS Publishing Company, Boston.

[48] STOLYAR, A. L. (1995). On the stability of multiclass queueing networks: A relaxed sufficient condition via limiting fluid processes. *Markov Process. Related Fields* **1** 491–512. MR1403094

[49] TSAPARAS, P. (1997). Stability in adversarial queueing theory. M.Sc. thesis, Univ. Toronto.
OPERATIONS RESEARCH CENTER AND
 SLOAN SCHOOL OF MANAGEMENT
MASSACHUSETTS INSTITUTE OF TECHNOLOGY
CAMBRIDGE, MASSACHUSETTS 02139
USA
E-MAIL: gamarnik@mit.edu

IBM T.J. WATSON RESEARCH CENTER
PO BOX 218
YORKTOWN HEIGHTS, NEW YORK 10598
USA
E-MAIL: dimdim@mit.edu